%% file: preprint.tex
\documentclass[twocolumn]{mpi2015-cscpreprint}

\providecommand\iffinal\iftrue 
\iffinal\PassOptionsToPackage{final}{changes}\fi
\usepackage[defaultcolor=orange]{changes}

\usepackage{manuscript}
\renewcommand\single{\makebox[\widthof{D}]{S}}
\renewcommand\double{D}
\usepackage{mycode}
\usepackage{orcidlink}
\declaretheorem{example}
\declaretheorem{remark}

\usepackage{microtype}

\usepackage{pgfplots}
\usepgfplotslibrary{groupplots}
\definecolor{colorbrewera}{HTML}{E66101}
\definecolor{colorbrewerb}{HTML}{FDB863}
\definecolor{colorbrewerc}{HTML}{B2ABD2}
\definecolor{colorbrewerd}{HTML}{5E3C99}
\pgfplotsset{
  compat=1.18,
  /tikz/implicit residual/.style = {
    mark = square*,
    mark size = 3pt,
    black!15,
  },
  every axis/.append style = {
    cycle list = {
      colorbrewera,mark=square*\\
      colorbrewerb,mark=*\\
      colorbrewerc,mark=pentagon*\\
      colorbrewerd,mark=triangle*\\
    },
    width = 0.275\linewidth,
    height = 9em,
    implicit residual/.style = {
      forget plot,
      /tikz/implicit residual,
    }
  },
  every axis legend/.append style = {
    /tikz/every even column/.append style = {
      column sep = 2ex,
    },
  },
  table/search path = {.,plots,plots-fallback},
}

\title{%
  Towards a mixed-precision ADI method\\ for Lyapunov equations
}

\author[$\ast\dagger$]{Jonas~Schulze\,\orcidlink{0000-0002-2086-7686}} 
\author[$\ast$]{Jens~Saak\,\orcidlink{0000-0001-5567-9637}} 

\affil[$\ast$]{Max Planck Institute for Dynamics of Complex Technical Systems, Magdeburg, Germany.}
\affil[$\dagger$]{Corresponding author.\ \email{jschulze@mpi-magdeburg.mpg.de}}

\abstract{\input{abstract}}
\novelty{%
  This work marks the first application of mixed-precision to the low-rank ADI.\@
  Our main focus is to keep the maximum memory required during the ADI as low as possible.
}

\begin{document}

\maketitle

\input{content}

\addcontentsline{toc}{section}{References}
\bibliographystyle{abbrvurl}
\bibliography{journals.bib,csc.bib,mor.bib,software.bib,literature.bib}

\end{document}

%% file: abstract.tex
We apply mixed-precision to the low-rank Lyapunov ADI (LR-ADI) by performing certain aspects of the algorithm in a lower working precision.
Namely, we accumulate the overall solution,
solve the linear systems comprising the ADI iteration,
and store the inner low-rank factors of the residuals in various combinations of IEEE~754 single and double precision.
We empirically test our implementation on Lyapunov equations arising from first- and second-order descriptor systems.
For the first-order examples,
accumulating the solution in single-precision yields an almost-as-small residual as for the double-precision solution.
For certain applications,
like computing the \Htwo{} norm of a descriptor system,
low- or mixed-precision variants of the ADI can be quite competitive.

%% file: content.tex

\section{Introduction}

We consider the numerical solution of the continuous-time \ac{ALE}
\begin{equation}%
\label{eq:intro}
  \Lyap(X) := A\transpose XE + E\transpose XA + GSG\transpose = 0
  ,
\end{equation}
with large and sparse coefficient matrix $A\in\C^{n\times n}$
and a low-rank constant term comprised of the factors $G\in\C^{n\times g}$ and $S\in\C^{g\times g}$,
where $g\ll n$ \added{and $S$ is Hermitian}.
This type of equation arises in, e.g., optimal control and model order reduction.
We refer to~\cite{morAnt05,BenS13,Sim16a} and the references therein for a more detailed introduction.

Given the low rank of the constant term in~\eqref{eq:intro},
the solution can, at least numerically, be well approximated by a low-rank factorization~\cite{Pen00,AntSZ02,Gra04,Sab07,TruV07}.
We choose the \replaced{Hermitian}{symmetric} indefinite factorization, \mbox{$X\approx \ZYZ\transpose$},
with tall and skinny $Z\in\C^{n\times z}$ and Hermitian $Y\in\C^{z\times z}$, $z\ll n$,
~\cite{BenLP08,LanMS15}.
For this kind of large-scale equations with low-rank solutions,
the most successful algorithms used recently are the Krylov subspace projection method,
\eg~\cite{Sim07,PalS18,KreLMetal21},
as well as the \ac{ADI} method;
see~\cite{SchS25} and references therein.
The latter will be the focus of this paper.

In recent years,
driven by the needs of machine learning,
hardware capable of computing in lower than IEEE~754~\cite{IEEE754} single precision has been becoming available to the general user.
Many algorithms from all fields of mathematics
have been found to yield just- or almost-as-good results
when using lower precision for parts of their computations,
while being faster or consuming less energy.
We refer to~\cite{AbdABetal21,HigM22} and the references therein for a more detailed introduction.
Motivated by~\cite{AliAGetal23},
as a first step towards lower than IEEE~754~\cite{IEEE754} single precision,
we seek to employ mixed-precision for aforementioned low-rank factorizations in the context of large-scale matrix equations.
For this paper, our motivation is to obtain same-quality solutions to \ac{ALE}~\eqref{eq:intro} while using less memory to store the solution.
A detailed examination of runtime and energy consumption is out of scope for this initial investigation,
and left for future work.

\subsection{Notation}

Throughout the paper,
$\Fnorm{\cdot}$ denotes the Frobenius norm.
Hermitian transposition of a matrix is denoted by $(\cdot)\transpose : \C^{m\times n} \to \C^{n\times m}$.
We consistently use the same notation for the transposition of a real matrix.
The identity matrix is denoted by $I_n\in\R^{n\times n}$,
where we may omit the subscript if dimensions are clear from the context.
Conceptually, we identify the zero matrix with a rank-zero decomposition $[\,]\,[\,]\,[\,]\transpose = 0\in\R^{n\times n}$.

We denote the machine epsilon of a data type by $\precision{}\in\R$ where $0<\precision{}\ll 1$.
Throughout the paper, we will only use IEEE~754~\cite{IEEE754} precisions single and double,
such that we may identify a data type with its machine epsilon.
Subscripts of \precision{} denote the matrices stored in the corresponding data type,
for example, \precTY{} refers the precision used to store matrices named $T$ and $Y$.\footnote{%
  The subscript does \emph{not} denote matrix multiplication.
}

\subsection{Low-rank arithmetic}%
\label{sec:arithmetic}

To compress a low-rank factorization $\ZYZ\transpose$,
we follow~\cite{LanMS15} and compute a thin QR-decomposition (pivoted Householder~\cite{Hou58,GolV13}) of $Z = Q\hat Z$ using precision~\precZ,
and truncate based on an eigendecomposition of $\hat Z Y_{} \hat Z\transpose$ using \precTY.\footnote{%
  $T$ will be denote an inner low-rank factor akin to $Y$,
  whose exact meaning will be introduced later.
}
To compute the norm of $\ZYZ\transpose$,
we follow, e.g.,~\cite{Pen99} and evaluate $\Fnorm{\smash{\hat Z Y_{} \hat Z\transpose}}$ using \precTY{}.
In this sense,
throughout the paper,
the working precision of certain operations will always follow the storage precision of the matrices involved.
We refer the reader to our implementation~\cite{DRE.jl} for further details.

We call $\ZYZ\transpose$ with $Z\in\R^{n\times z}$ an order-$z$ approximation.
Note that $z$ is only an upper bound to $\rank(\ZYZ\transpose)$.
After low-rank compression as described above,
the order and rank of a low-rank factorization coincide.

\section{Mixed-precision low-rank factorization}

Let $\sizeof(\cdot)$ denote the number of Bytes needed to store a certain object.
Furthermore, let $\ZYZ\transpose$ be a proper low-rank factorization,
i.e., $Z\in\R^{n\times z}$ with $z\ll n$.
Assuming uniform precision and dense storage for~$Z$ and~$Y$,
it holds
\begin{equation}
  \sizeof(Z) \gg \sizeof(Y)
  .
\end{equation}
A natural choice to reduce the total storage required
is to permit a coarser precision for the factor~$Z$.

Now let $\ZYZ\transpose$ be an approximation of the true solution~$X$ to \ac{ALE}~\eqref{eq:intro} obtained via the \ac{ADI} method,
and let $\RTR\transpose$ be the corresponding residual, i.e.,
\begin{math}
  \Lyap(\ZYZ\transpose) = \RTR\transpose
  .
\end{math}
Recall from~\cite[Algorithm~3.1]{LanMS15}
that the inner solution factor~$Y$ can be represented as a Kronecker product,
\begin{equation}%
\label{eq:kronecker}
  Y = -2\operatorname{diag}(\Re(\alpha_0), \ldots, \Re(\alpha_k)) \otimes T
  ,
\end{equation}
\added{%
while the outer factor~$Z$ decomposes into equal-sized blocks,
\begin{math}
  Z = [V_0 | \ldots | V_k].
\end{math}
}
Thus, assuming uniform precision, $z \leq n$,
and dense storage for all matrices but $Y$,
we conclude
\begin{equation}
\begin{alignedat}{2}
  \sizeof(Z)
  &\geq \sizeof(V_{\added{*}}) &&= \sizeof(R) \\
  &\geq \sizeof(T) &&= \mathcal O(\sizeof(Y))
  ,
\end{alignedat}
\end{equation}
where $\mathcal O(\cdot)$ denotes Landau notation.
Similarly,
a natural choice to minimize the overall storage required
is to permit three precisions
\begin{equation}
  \precZ \geq \precVR \geq \precTY
  .
\end{equation}

\Autoref{alg:adi} shows a simplified summary of the resulting mixed-precision \ac{ADI}.
For brevity, we omit the handling of complex shifts~\cite{BenKS13}.
To solve the linear systems in \autoref{alg:adi:linsolve},
we use a column-wise applied GMRES~\cite{SaaS86,Krylov.jl}
with ILU(0) preconditioner~\cite{MeiV77,ILUZero.jl}. 

\begin{remark}%
\label{rem:GMRES}
  A current limitation of our implementation of \autoref{alg:adi},
  see DifferentialRiccatiEquations.jl~\cite{DRE.jl},
  is that the matrices~$E$ and~$A$ as well as the shifts~$\alpha_k$
  must all be converted to~\precVR{} beforehand.
  This conversion is necessary to ensure type stability,
  and thus reduce the number of allocations inside Krylov.jl~\cite{Krylov.jl},
  the package that implements GMRES\@.
\end{remark}

\input{adi_algorithm}

\section{Numerical experiments}

In this section,
we investigate the following questions:
\begin{enumerate}
  \item
    Can we approximate the solutions with maximum attainable accuracy?
  \item
    How will convergence be affected by using lower working precisions?
\end{enumerate}
We aim to answer these questions separately
for typical problems arising in control theory
(\autoref{sec:num:gramian}),
and in a more general setting
(\autoref{sec:num:random-rhs} and \autoref{sec:num:known-X}).
Our primary motivation is to obtain solutions of the same quality as if using double precision.

An \ac{ALE}~\eqref{eq:intro} can be derived from a linear state-space system
\begin{equation}%
  \label{eq:EABC}
  \Sigma:\left\{
  \quad
  \begin{aligned}
    E\dot x(t) &= Ax(t) + Bu(t), \\
    y(t) &= Cx(t),
  \end{aligned}
  \right.
\end{equation}
with states~$x:\R\to\R^n$, inputs~$u:\R\to\R^m$, and outputs~$y:\R\to\R^q$.
For the remainder of this section,
we consider data from \added{the} following underlying dynamical systems:

\begin{example}[Steel Profile~\cite{morwiki_steel,BenS05b,SaaB21}]
  This system describes a semi-discretized heat transfer problem for optimal cooling of steel profiles.
  The matrices~$E$ and~$A$ are both symmetric,
  and the dimensions are
  \begin{equation*}
    n=\num{5177},
    \quad
    m=7,
    \quad
    q=6.
  \end{equation*}
  The data set has been downloaded via
  \lstinline{FenicsRail(5177)} 
  using MORWiki.jl~\cite{MORWiki.jl}.
\end{example}

\begin{example}[Chip~\cite{morwiki_convection,morMooRGetal04}]
  This system describes cooling of a computer chip.
  The system matrix~$E$ is symmetric, while matrix~$A$ is non-symmetric.
  The dimensions are
  \begin{equation*}
    n=\num{20082},
    \quad
    m=1,
    \quad
    q=5.
  \end{equation*}
  The data set has been downloaded via
  \lstinline{Chip(0.0)} 
  using MORWiki.jl~\cite{MORWiki.jl}.
\end{example}

\begin{example}[Triple Chain~\cite{TruV09}]%
\label{ex:triple chain}
  This system describes
  three parallel equally sized chains of masses, springs, and dampers,
  that are coupled by one big connector mass.
  The system is modeled by a second-order formulation
  \begin{equation*}
    \left\{
    \begin{aligned}
      \hat M\ddot x(t) + \hat E\dot x(t) + \hat Kx(t) &= \hat Bu(t), \\
      y(t) &= \hat Cx(t),
    \end{aligned}
    \right.
  \end{equation*}
  where $\hat M, \hat E$, and $\hat K\in\R^{\hat n\times\hat n}$ describe mass, damping,\footnote{%
    In general, the matrix~$\hat E$ describes damping and gyroscopic effects.
    However, there are no gyroscopic effects here.
  }
  and stiffness, respectively.
  We use Rayleigh damping, $\hat E=\alpha \hat M + \beta \hat K$,
  with $\alpha=\beta=0.1$ for a viscosity of $v=5$.
  The other system parameters as described in~\cite[Example~2]{TruV09} are given by
  \begin{equation*}
    \begin{aligned}
      k_0 &= 50, &
      k_1 &= 10, &
      k_2 &= 20, &
      k_3 &= 1, \\
      m_0 &= 10, &
      m_1 &= 1, &
      m_2 &= 2, &
      m_3 &= 3.
    \end{aligned}
  \end{equation*}
  An equivalent first-order formulation~\eqref{eq:EABC} of dimension $n=2\hat n$
  has been derived in a strictly-dissipative way;
  see~\cite{morPan14} and references therein.
  The resulting system matrix $E$ is symmetric,
  $A\in\R^{n\times n}$ is non-symmetric,
  and the dimensions are
  \begin{equation*}
    n=\num{602},
    \quad
    m=1,
    \quad
    q=1.
  \end{equation*}
  This system may seem small, however, it will prove to be difficult enough to solve.
\end{example}

To isolate the storage/working precision as the only difference between the algorithms,
we use the same pre-computed shifts~\cite{Pen99} for all variants of our \ac{ADI};
obtained via
\lstinline{Heuristic(20, 40, 40)} 
using DifferentialRiccatiEquations.jl~\cite{DRE.jl};
and request a relative tolerance of \replaced{$\tau=10^{-10}$}{$\tau=10^{-8}$} regardless of working precision.
For the Chip and the Steel Profile examples,
we allow a maximum of 50 \ac{ADI} iterations,
whereas for the Triple Chain we allow 200 iterations.\footnote{%
  None of our experiments will exhaust this limit.
  However,
  the maximum number of iterations allowed is part of the file names we use to store the results.
  To avoid excessive computation times when trying to reproduce the results,
  we select a limit that is only slightly larger than the actually needed number of iterations.
}
We request the inner GMRES to compute solutions to a relative tolerance of $100\precVR$.

\begin{remark}%
\label{rem:GMRES:precision}
  The explicit Lyapunov residuals are quite sensitive to the accuracy of the increments (\autoref{alg:adi:linsolve} in \autoref{alg:adi}),
  whereas the implicit residuals are not.
  Even for the uniform double-precision ADI,
  using the default relative tolerance of $\sqrt{\precVR}$ for the GMRES as implemented in Krylov.jl~\cite{Krylov.jl},
  causes the explicit residuals to diverge from the implicit ones.
  At some point of the iteration~$k\in\N$,
  the explicit residuals stop decreasing.
  If we instead use a GMRES tolerance of $100\precVR$,
  we only observe this divergence for much smaller ADI tolerances~$\tau$.
  \added{%
  The uniform double-precision ADI, for example,
  will only show a stagnating explicit residual in \autoref{sec:num:gramian}.
  In fact, if we had used a direct inner solver for the first couple of \ac{ADI} iterations,
  the explicit residual would have converged;
  see the Pluto.jl~\cite{Pluto.jl} notebook supplementing our codes.
  More rigorous tolerances for the inner solver can be chosen based on~\cite{KueF20}.%
  }%
\end{remark}

The following subsections show three numerical experiments.
First, we solve a very specific application of the \ac{ALE},
where one is typically not interested in the solution itself
but in certain quantities that can be derived from the solution without actually assembling it,
here demonstrated by the computation of the \Htwo{} system norm.
As the exact solution is not known,
we compare our iterative \ac{ADI} algorithm to the direct Bartels-Stewart algorithm.
Second, we solve a more general type of \ac{ALE} with a random constant term of varying rank.
Again, the exact solution is not known,
but given the findings from the first experiment and to save some electricity,
we will not apply the Bartels-Stewart algorithm.
Third, we solve an \ac{ALE} with a known solution of varying rank.
In the above,
the experiments were mainly concerned with the residual,
whereas here we also investigate the solution error.

Our expectation is that the four instances of ADI(\precZ, \precVR, \precTY) as presented by \autoref{alg:adi} ordered by increasing use of IEEE single precision;
ADI(D, D, D),
ADI(S, D, D),
ADI(S, S, D), and
ADI(S, S, S);
are ordered by decreasing solution quality and decreasing runtime.
Our hypotheses are that
\begin{enumerate}
  \item\label{list:hypotheses}
    Overall convergence will be the same,
  \item
    \mbox{ADI(S, D, D)} will yield solutions of the same quality
    as the uniform double-precision variant while being faster,
    and that
  \item
    \mbox{ADI(S, S, D)} will yield solutions of better quality
    than the uniform single-precision variant
    at essentially the same runtime.
\end{enumerate}

We assess the quality of a low-rank \ac{ADI} iterate~$Z_k Y_k Z_k\transpose$
by means of the normalized \emph{implicit} (internally updated) and \emph{explicit residuals} (inserted in~\eqref{eq:intro}),
defined by
\begin{equation}%
\label{eq:residuals}
  \resimpl{R_k T R_k\transpose}
  \quad\text{and}\quad
  \resexpl{Z_k Y_k Z_k\transpose}
  ,
\end{equation}
where the residual factors~$R_k$ and~$T$ are given by \autoref{alg:adi},
and all Frobenius norms are evaluated as described in \autoref{sec:arithmetic}.
Recall from, e.g.,~\cite{SchS25} that
\begin{math}
  \Lyap(\ZYZ\transpose) = \hat R \hat T \hat R\transpose
\end{math}
with
\begin{equation}
  \hat R = \begin{bmatrix}
    G & EZ & AZ
  \end{bmatrix}
  \quad\text{and}\quad
  \hat T = \begin{bmatrix}
    S \\
    && Y \\
    & Y
  \end{bmatrix}
  .
\end{equation}

\replaced{%
All experiments were performed on Ubuntu 24.04.3 with kernel 6.14.0-37-generic, 
running on an AMD EPYC\textsuperscript{\tiny TM} 9554 processor with 64~cores.
}{%
Most experiments were performed on an Intel\textsuperscript{\tiny\circledR} Core\textsuperscript{\tiny TM} i5-12600K with 6~P-cores and 4~E-cores. 
The only exception is the Bartels-Stewart method applied to the Chip example,
for which we used an AMD EPYC\textsuperscript{\tiny TM}~9554 processor with 64~cores.
The increased computational power provided by the AMD processor was necessary to ensure the computations were completed in a reasonable amount of time.
}

\subsection{Observability Gramian and \texorpdfstring{\Htwo}{H2} norm}%
\label{sec:num:gramian}

The observability Gramian~$E\transpose QE\in\R^{n\times n}$ can be computed by means of the solution~$Q$ to
\begin{equation}
  \label{eq:gramian}
  A\transpose Q E + E\transpose Q A = -C\transpose C,
\end{equation}
and has high practical relevance to the fields of systems and control theory as well as model order reduction.
In terms of \ac{ALE}~\eqref{eq:intro},
we set the right-hand side to $G=C\transpose\in\R^{n\times q}$ and $S=I_q$.
Another, derived interesting property is the \Htwo{} norm of system~\eqref{eq:EABC};
see, e.g.,~\cite{morAnt05};
and can be computed via
\begin{equation}
  \Hnorm{\Sigma}^2 =
  \operatorname{trace}(B\transpose Q B)
  .
\end{equation}

\Autoref{tab:gramian}
compares the solutions obtained via
\autoref{alg:adi}
to a (densely computed) Bartels-Stewart reference solution implemented in MatrixEquations.jl~\cite{BarS72,ME.jl}.
Observe that \added{almost} all implicit residuals~\eqref{eq:residuals}
reach the desired tolerance of \replaced{$\tau=10^{-10}$}{$\tau=10^{-8}$}.
However, the explicit residuals~\eqref{eq:residuals}
(after casting the low-rank factors of $Q$ to IEEE double precision) reveals that
only the uniform double-precision ADI and Bartels-Stewart algorithms
yield an accurate solution to equation~\eqref{eq:gramian}.
As soon as IEEE single precision is involved;
with the exception of \mbox{ADI(S, D, D)} applied to the \deleted{Steel Profile or} Chip example;
the explicit residual is multiple orders of magnitude larger than the implicit one.
Nevertheless, all derived \Htwo{} norms for the Steel Profile essentially coincide,
while for the Chip and Triple Chain examples,
the \Htwo{} norms of the first two ADI variants are basically the same as for the double-precision Bartels-Stewart algorithm.
Meanwhile, the \Htwo{} norms of the last two ADI variants differs from the reference value by less than
\qty{3.6}{\percent} and \qty{0.1}{\percent} \added{for the Chip and Triple Chain examples, respectively}.
In comparison, the single-precision Bartels-Stewart algorithm yields \Htwo{} norms for the Chip and Triple Chain examples
that differ from their double-precision values by about \qty{55}{\percent} and \qty{3.9}{\percent}, respectively,
while their corresponding explicit residuals (evaluated in double precision) are \numrange{9}{10}~orders larger;
much larger even than for the uniformly single-precision ADI\@.

Note that, in general, the reconstruction errors have about the same magnitude as the explicit residuals.\footnote{%
  The exception is given by \mbox{ADI(S, $\ast$, $\ast$)} applied to the Triple Chain,
  where the error is about \numrange{4}{6}~orders smaller, i.e.,~better than expected.
}
Also,
while the runtime of both \mbox{ADI($\ast$, D, D)} and both \mbox{ADI(S, S, $\ast$)} are the same,
the Bartels-Stewart method takes orders of magnitude longer.
For this reason, we will not employ Bartels-Stewart in later numerical experiments.

\Autoref{fig:num:gramian} shows the implicit and explicit residuals over the course of the ADI iteration.
Observe that for all examples,
\added{almost}
all the implicit residuals coincide,
while the explicit residual trajectories agree with the implicit ones only up to a certain iteration;
see \autoref{rem:GMRES:precision}.
Afterwards, the explicit residual stagnates.
\added{%
The only exception is given by the Chip example,
for which even the implicit residuals of the low-precision \mbox{ADI(S, S, $\ast$)} stagnate.%
}\footnote{\added{%
  Due to stagnation slightly above the target tolerance of $\tau=10^{-10}$,
  the method proceeds until exhausting the maximum 50~iterations,
  or until the increment vanishes, that is, $V_k = 0$.
}}\ %
\added{%
Note that even the explicit \mbox{ADI(D, D, D)} residuals stagnate;
except for the Steel Profile.
}%
Given that the residuals of all the variants of the \ac{ADI} behave this way,
to reduce the energy consumption of the remaining numerical experiments in this section,
we only focus on the explicit residuals~\eqref{eq:residuals}
for the last iteration~$k\in\N$.

\begin{table*}
  \footnotesize
  \centering
  \caption{%
    Numerical results of computing the Observability Gramian~\eqref{eq:gramian};
    see \autoref{sec:num:gramian}.
    The symbols~\single{} and~\double{} refer to IEEE single and double precision, respectively.
    The remaining columns show the orders of the
    final ADI residual~$\RTR\transpose$ and
    solution~$\ZYZ\transpose$ representing the Gramian~$Q$;
    the corresponding normalized implicit and explicit residual;
    the normalized error versus double-precision Bartels-Stewart solution~$\hat Q$;
    the system's \Htwo{} norm;
    the total number of \ac{ADI} iterations;
    and the wall-clock time of the solver.
  }%
  \label{tab:gramian}
  \sisetup{%
    exponent-mode=scientific,
    round-mode=places,
    round-precision=2,
    tight-spacing=true,
  }
  \begin{tabular}{%
      @{}
      ll 
      S[table-format=1,round-mode=none,exponent-mode=fixed,fixed-exponent=0] 
      S[table-format=3,round-mode=none,exponent-mode=fixed,fixed-exponent=0] 
      S[table-format=1.2e-1] 
      S[table-format=1.2e-2] 
      S[table-format=1.2e-2] 
      S[table-format=1.4e-1,round-precision=4] 
      S[table-format=3,round-mode=none,exponent-mode=fixed,fixed-exponent=0] 
      S[table-format=3.2{\,\unit{\minute}},table-align-text-after=false,exponent-mode=fixed,fixed-exponent=0] 
      @{}
    }
    \toprule
    &&
    \multicolumn{2}{c}{Order} &
    \multicolumn{2}{c}{Residual} &
    \\
    \cmidrule(rl){3-4}
    \cmidrule(rl){5-6}
    Example &
    Method &
    {$\RTR\mathrlap{\transpose}$} &
    {$\ZYZ\mathrlap{\transpose}$} &
    {\reslabelimpl{\smash{\RTR\transpose}}} &
    {\reslabelexpl{Q}} &
    {$\frac{\Fnorm{\smash{Q-\hat Q}}}{\Fnorm{\smash{\hat Q}}}$} &
    {$\Hnorm{\Sigma}$} &
    {\#Iter.} &
    {Time}
    \\
    \midrule
    Steel Profile
    & ADI(\double, \double, \double) & 6 & 222 & 6.629068956382094e-11 & 6.629844343203268e-11 & 6.433346450341527e-11 & 1.0407304561779315e-5 & 37 & 3.2438127420000002\,\si{\second} \\
    & ADI(\single, \double, \double) & 6 & 222 & 6.629068956382094e-11 & 5.682703754236161e-8 & 2.7137989406084586e-8 & 1.0407304562321484e-5 & 37 & 3.2292034010000004\,\si{\second} \\
    & ADI(\single, \single, \double) & 6 & 222 & 6.628611825887971e-11 & 1.772565646651992e-5 & 4.385987564254357e-5 & 1.0407297211817498e-5 & 37 & 819.775594\,\si{\milli\second} \\
    & ADI(\single, \single, \single) & 6 & 222 & 6.628611939204665e-11 & 1.772565646651992e-5 & 4.386401935925323e-5 & 1.0407297211817498e-5 & 37 & 818.108137\,\si{\milli\second} \\
    & Bartels-Stewart(\double) & {--} & {--} & {--} & 2.6109089516588942e-14 & 0.0 & 1.0407304561934519e-5 & {--} & 14.1650637971\,\si{\minute} \\
    & Bartels-Stewart(\single) & {--} & {--} & {--} & 1.6330534429131586e-5 & 0.00014840854142706943 & 1.040795813676507e-5 & {--} & 12.5453610764\,\si{\minute} \\
    \midrule
    Chip
    & ADI(\double, \double, \double) & 5 & 160 & 2.1894736906768104e-11 & 8.079832185209008e-9 & 8.076877356661235e-11 & 323.53229499195504 & 32 & 8.699200737\,\si{\second} \\
    & ADI(\single, \double, \double) & 5 & 160 & 2.1894736906768104e-11 & 5.387676828186326e-8 & 2.498923338585434e-8 & 323.5322949190987 & 32 & 8.683178438\,\si{\second} \\
    & ADI(\single, \single, \double) & 5 & 210 & 3.291549594800563e-10 & 0.0003133482534714006 & 1.6462038095665466e-5 & 312.0007078208899 & 43 & 2.752598795\,\si{\second} \\
    & ADI(\single, \single, \single) & 5 & 255 & 3.2915495824565385e-10 & 0.00031334825347139727 & 1.646436985932219e-5 & 312.0007078208899 & 51 & 2.9118611420000002\,\si{\second} \\
    & Bartels-Stewart(\double) & {--} & {--} & {--} & 1.55597342865277e-10 & 0.0 & 323.5322987971623 & {--} & 18.347489828919723\,\si{\hour} \\
    & Bartels-Stewart(\single) & {--} & {--} & {--} & 1.3002943267202507 & 0.6489317432311847 & 502.40346157182967 & {--} & 15.048747532523333\,\si{\hour} \\
    \midrule
    Triple Chain
    & ADI(\double, \double, \double) & 1 & 201 & 1.4036855224016042e-11 & 2.342312034564651e-10 & 7.355222186704727e-10 & 5065.677210238942 & 201 & 22.584192\,\si{\milli\second} \\
    & ADI(\single, \double, \double) & 1 & 201 & 1.4036855224012617e-11 & 0.03823350152511256 & 3.2949753055328256e-8 & 5065.677202396358 & 201 & 22.898161\,\si{\milli\second} \\
    & ADI(\single, \single, \double) & 1 & 191 & 2.4011224793543774e-10 & 0.1239615154860279 & 3.449146788469654e-5 & 5065.54483550407 & 193 & 20.07135\,\si{\milli\second} \\
    & ADI(\single, \single, \single) & 1 & 191 & 2.4011224311901143e-10 & 0.12396151548602805 & 3.4555597854692954e-5 & 5065.544835504076 & 193 & 19.880494\,\si{\milli\second} \\
    & Bartels-Stewart(\double) & {--} & {--} & {--} & 1.1132051235964606e-7 & 0.0 & 5065.677209354674 & {--} & 1.3988422539999998\,\si{\second} \\
    & Bartels-Stewart(\single) & {--} & {--} & {--} & 44.49836054818302 & 0.20051156811472673 & 5292.617937036107 & {--} & 805.025474\,\si{\milli\second} \\
    \bottomrule
  \end{tabular}
\end{table*}

\begin{figure*}
  \centering
  \begin{tikzpicture}
    \pgfplotstableread[col sep=comma]{Chip_0.0_/gramian/ADI-DDD_compression=false_compression_interval=10_ignore_initial_guess=false_maxiters=50_reltol=1e-10_.csv}\adiDDDchip
    \pgfplotstableread[col sep=comma]{Chip_0.0_/gramian/ADI-SDD_compression=false_compression_interval=10_ignore_initial_guess=false_maxiters=50_reltol=1e-10_.csv}\adiSDDchip
    \pgfplotstableread[col sep=comma]{Chip_0.0_/gramian/ADI-SSD_compression=false_compression_interval=10_ignore_initial_guess=false_maxiters=50_reltol=1e-10_.csv}\adiSSDchip
    \pgfplotstableread[col sep=comma]{Chip_0.0_/gramian/ADI-SSS_compression=false_compression_interval=10_ignore_initial_guess=false_maxiters=50_reltol=1e-10_.csv}\adiSSSchip
    \pgfplotstableread[col sep=comma]{FenicsRail_5177_/gramian/ADI-DDD_compression=false_compression_interval=10_ignore_initial_guess=false_maxiters=50_reltol=1e-10_.csv}\adiDDDrail
    \pgfplotstableread[col sep=comma]{FenicsRail_5177_/gramian/ADI-SDD_compression=false_compression_interval=10_ignore_initial_guess=false_maxiters=50_reltol=1e-10_.csv}\adiSDDrail
    \pgfplotstableread[col sep=comma]{FenicsRail_5177_/gramian/ADI-SSD_compression=false_compression_interval=10_ignore_initial_guess=false_maxiters=50_reltol=1e-10_.csv}\adiSSDrail
    \pgfplotstableread[col sep=comma]{FenicsRail_5177_/gramian/ADI-SSS_compression=false_compression_interval=10_ignore_initial_guess=false_maxiters=50_reltol=1e-10_.csv}\adiSSSrail
    \pgfplotstableread[col sep=comma]{TripleChain_100_/gramian/ADI-DDD_compression=false_compression_interval=10_ignore_initial_guess=false_maxiters=200_reltol=1e-10_.csv}\adiDDDmsd
    \pgfplotstableread[col sep=comma]{TripleChain_100_/gramian/ADI-SDD_compression=false_compression_interval=10_ignore_initial_guess=false_maxiters=200_reltol=1e-10_.csv}\adiSDDmsd
    \pgfplotstableread[col sep=comma]{TripleChain_100_/gramian/ADI-SSD_compression=false_compression_interval=10_ignore_initial_guess=false_maxiters=200_reltol=1e-10_.csv}\adiSSDmsd
    \pgfplotstableread[col sep=comma]{TripleChain_100_/gramian/ADI-SSS_compression=false_compression_interval=10_ignore_initial_guess=false_maxiters=200_reltol=1e-10_.csv}\adiSSSmsd
    \begin{groupplot}[
      group style = {
        columns = 3,
        rows = 2,
        horizontal sep = 2ex,
        vertical sep = 2ex,
        x descriptions at = edge bottom,
        y descriptions at = edge left,
      },
      scale only axis,
      ymode = log,
      ytickten = {0, -5, -10},
      ymax = 20,
      ymin = 1e-11,
      enlargelimits = true,
      /tikz/only marks,
    ]
    \nextgroupplot[
      title = Steel Profile,
      legend to name = fig:num:gramian:legend,
      legend columns = -1,
      legend entries = {%
        {ADI(\double, \double, \double)},
        {ADI(\single, \double, \double)},
        {ADI(\single, \single, \double)},
        {ADI(\single, \single, \single)},
      },
      ylabel = \reslabelimpl{R_k T R_k},
    ]
      \addplot table [x="iter",y="res_impl"] {\adiDDDrail};
      \addplot table [x="iter",y="res_impl"] {\adiSDDrail};
      \addplot table [x="iter",y="res_impl"] {\adiSSDrail};
      \addplot table [x="iter",y="res_impl"] {\adiSSSrail};
    \nextgroupplot[title=Chip, xmax=32]
      \addplot table [x="iter",y="res_impl"] {\adiDDDchip};
      \addplot table [x="iter",y="res_impl"] {\adiSDDchip};
      \addplot table [x="iter",y="res_impl"] {\adiSSDchip};
      \addplot table [x="iter",y="res_impl"] {\adiSSSchip};
    \nextgroupplot[title=Triple Chain, xmax=201]
      \addplot table [x="iter",y="res_impl"] {\adiDDDmsd};
      \addplot table [x="iter",y="res_impl"] {\adiSDDmsd};
      \addplot table [x="iter",y="res_impl"] {\adiSSDmsd};
      \addplot table [x="iter",y="res_impl"] {\adiSSSmsd};
    \nextgroupplot[ 
      ylabel = \reslabelexpl{X_k},
    ]
      \addplot[implicit residual] table [x="iter",y="res_impl"] {\adiDDDrail};
      \addplot table [x="iter",y="res_expl"] {\adiDDDrail};
      \addplot table [x="iter",y="res_expl"] {\adiSDDrail};
      \addplot table [x="iter",y="res_expl"] {\adiSSDrail};
      \addplot table [x="iter",y="res_expl"] {\adiSSSrail};
    \nextgroupplot[ 
      xlabel = Iteration $k$,
      xmax = 32,
    ]
      \addplot[implicit residual] table [x="iter",y="res_impl"] {\adiDDDchip};
      \addplot table [x="iter",y="res_expl"] {\adiDDDchip};
      \addplot table [x="iter",y="res_expl"] {\adiSDDchip};
      \addplot table [x="iter",y="res_expl"] {\adiSSDchip};
      \addplot table [x="iter",y="res_expl"] {\adiSSSchip};
    \nextgroupplot[xmax=201] 
      \addplot[implicit residual] table [x="iter",y="res_impl"] {\adiDDDmsd};
      \addplot table [x="iter",y="res_expl"] {\adiDDDmsd};
      \addplot table [x="iter",y="res_expl"] {\adiSDDmsd};
      \addplot table [x="iter",y="res_expl"] {\adiSSDmsd};
      \addplot table [x="iter",y="res_expl"] {\adiSSSmsd};
    \end{groupplot}
  \end{tikzpicture}

  \medskip
  \ref*{fig:num:gramian:legend} 
  \caption{%
    Implicit (top) and explicit residuals (bottom)
    of computing the Observability Gramian~\eqref{eq:gramian} over iteration index~$k\in\N$ for
    ADI(\precZ, \precVR, \precTY) 
    as shown in \autoref{alg:adi};
    see \autoref{sec:num:gramian}.
    The bottom row shows implicit residuals in the background~(%
    \protect\tikz[/pgfplots/every crossref picture] \protect\draw[implicit residual] plot coordinates { (0,0) };).
    The symbols S and D refer to IEEE single and double precision, respectively.
  }%
  \label{fig:num:gramian}
\end{figure*}

\begin{remark}%
\label{rem:compression}
  The solutions reported in this subsection are all \emph{not} compressed.
  We observed low-rank compression following~\cite{LanMS15}
  applied to any $\ZYZ\transpose$ involving a single-precision pivoted Householder QR in~$Z$
  to increase the residual~$\Fnorm{\Lyap(\ZYZ\transpose)}$ by several orders of magnitude,
  even for small compression tolerances.
  This effect is much less pronounced if using other orthonormalization methods,
  for example, a shifted Cholesky QR algorithm.
  The residuals do not increase when casting the (previously truncated) data to double-precision before performing compression.
\end{remark}

\subsection{Random constant term}%
\label{sec:num:random-rhs}

In this experiment,
we \deleted{set} select a random constant term of \ac{ALE}~\eqref{eq:intro} for a given rank~$g\in\N$, $1\leq g \leq 50$.
In contrast to the Gramian described in \autoref{sec:num:gramian},
here the factor $S\in\R^{g\times g}$ is dense and indefinite.

\Autoref{fig:num:random-rhs} shows the implicit and explicit residuals for varying ranks~$g$.
Again, all implicit residuals suggest convergence to the desired tolerance of \replaced{$\tau=10^{-10}$}{$\tau=10^{-8}$},
while the explicit residuals reveal different behavior.
For any fixed variant of the ADI,
we observe a uniform quality of solutions over the rank~$g$.
Note that, as expected, the double-precision ADI produces the best results.
Accumulating the overall solution factor~$Z$ in single-precision
increases the residual norm by about \replaced{4, 4, and 6}{2, 2, and 4}~orders of magnitude for the Steel Profile, Chip, and Triple Chain, respectively.
Additionally solving the linear systems comprising the increments (\autoref{alg:adi:linsolve} of \autoref{alg:adi}) in single-precision
increases the residual norm by another 2~orders, 3~orders, and a factor of about~\replaced{3}{4}, respectively.
Truncating the common increment- and residual-factor~$T$ to single precision does not further increase the residual.
In terms of runtime, again,
both \mbox{ADI($\ast$, D, D)} and both \mbox{ADI(S, S, $\ast$)} are the same.

\begin{figure*}
  \centering
  \begin{tikzpicture}
    \pgfplotstableread[col sep=comma]{Chip_0.0_/random-rhs/ADI-DDD_compression=false_compression_interval=10_ignore_initial_guess=false_maxiters=50_reltol=1e-10_.csv}\adiDDDchip
    \pgfplotstableread[col sep=comma]{Chip_0.0_/random-rhs/ADI-SDD_compression=false_compression_interval=10_ignore_initial_guess=false_maxiters=50_reltol=1e-10_.csv}\adiSDDchip
    \pgfplotstableread[col sep=comma]{Chip_0.0_/random-rhs/ADI-SSD_compression=false_compression_interval=10_ignore_initial_guess=false_maxiters=50_reltol=1e-10_.csv}\adiSSDchip
    \pgfplotstableread[col sep=comma]{Chip_0.0_/random-rhs/ADI-SSS_compression=false_compression_interval=10_ignore_initial_guess=false_maxiters=50_reltol=1e-10_.csv}\adiSSSchip
    \pgfplotstableread[col sep=comma]{FenicsRail_5177_/random-rhs/ADI-DDD_compression=false_compression_interval=10_ignore_initial_guess=false_maxiters=50_reltol=1e-10_.csv}\adiDDDrail
    \pgfplotstableread[col sep=comma]{FenicsRail_5177_/random-rhs/ADI-SDD_compression=false_compression_interval=10_ignore_initial_guess=false_maxiters=50_reltol=1e-10_.csv}\adiSDDrail
    \pgfplotstableread[col sep=comma]{FenicsRail_5177_/random-rhs/ADI-SSD_compression=false_compression_interval=10_ignore_initial_guess=false_maxiters=50_reltol=1e-10_.csv}\adiSSDrail
    \pgfplotstableread[col sep=comma]{FenicsRail_5177_/random-rhs/ADI-SSS_compression=false_compression_interval=10_ignore_initial_guess=false_maxiters=50_reltol=1e-10_.csv}\adiSSSrail
    \pgfplotstableread[col sep=comma]{TripleChain_100_/random-rhs/ADI-DDD_compression=false_compression_interval=10_ignore_initial_guess=false_maxiters=200_reltol=1e-10_.csv}\adiDDDmsd
    \pgfplotstableread[col sep=comma]{TripleChain_100_/random-rhs/ADI-SDD_compression=false_compression_interval=10_ignore_initial_guess=false_maxiters=200_reltol=1e-10_.csv}\adiSDDmsd
    \pgfplotstableread[col sep=comma]{TripleChain_100_/random-rhs/ADI-SSD_compression=false_compression_interval=10_ignore_initial_guess=false_maxiters=200_reltol=1e-10_.csv}\adiSSDmsd
    \pgfplotstableread[col sep=comma]{TripleChain_100_/random-rhs/ADI-SSS_compression=false_compression_interval=10_ignore_initial_guess=false_maxiters=200_reltol=1e-10_.csv}\adiSSSmsd
    \begin{groupplot}[
      group style = {
        columns = 3,
        rows = 3,
        horizontal sep = 2ex,
        vertical sep = 2ex,
        x descriptions at = edge bottom,
        y descriptions at = edge left,
      },
      scale only axis,
      xmode = log,
      ymode = log,
      ymin = 1e-11,
      ymax = 1e-9,
      enlarge y limits = true,
      /tikz/only marks,
    ]
    \nextgroupplot[
      title = Steel Profile,
      legend to name = fig:num:random-rhs:legend,
      legend columns = -1,
      legend entries = {%
        {ADI(\double, \double, \double)},
        {ADI(\single, \double, \double)},
        {ADI(\single, \single, \double)},
        {ADI(\single, \single, \single)},
      },
      ylabel = \reslabelimpl{\RTR\transpose},
    ]
      \addplot table [x="rank",y="res_impl"] {\adiDDDrail};
      \addplot table [x="rank",y="res_impl"] {\adiSDDrail};
      \addplot table [x="rank",y="res_impl"] {\adiSSDrail};
      \addplot table [x="rank",y="res_impl"] {\adiSSSrail};
    \nextgroupplot[title = Chip]
      \addplot table [x="rank",y="res_impl"] {\adiDDDchip};
      \addplot table [x="rank",y="res_impl"] {\adiSDDchip};
      \addplot table [x="rank",y="res_impl"] {\adiSSDchip};
      \addplot table [x="rank",y="res_impl"] {\adiSSSchip};
    \nextgroupplot[title = Triple Chain]
      \addplot table [x="rank",y="res_impl"] {\adiDDDmsd};
      \addplot table [x="rank",y="res_impl"] {\adiSDDmsd};
      \addplot table [x="rank",y="res_impl"] {\adiSSDmsd};
      \addplot table [x="rank",y="res_impl"] {\adiSSSmsd};
    \nextgroupplot[ 
      ylabel = \reslabelexpl{X},
      ymin = 1e-11, ymax = 1e-2, ytickten = {-10, -8, ..., -2},
    ]
      \addplot[implicit residual] table [x="rank",y="res_impl"] {\adiDDDrail};
      \addplot table [x="rank",y="res_expl"] {\adiDDDrail};
      \addplot table [x="rank",y="res_expl"] {\adiSDDrail};
      \addplot table [x="rank",y="res_expl"] {\adiSSDrail};
      \addplot table [x="rank",y="res_expl"] {\adiSSSrail};
    \nextgroupplot[ 
      ymin = 1e-11, ymax = 1e-2, ytickten = {-10, -8, ..., -2},
    ]
      \addplot[implicit residual] table [x="rank",y="res_impl"] {\adiDDDchip};
      \addplot table [x="rank",y="res_expl"] {\adiDDDchip};
      \addplot table [x="rank",y="res_expl"] {\adiSDDchip};
      \addplot table [x="rank",y="res_expl"] {\adiSSDchip};
      \addplot table [x="rank",y="res_expl"] {\adiSSSchip};
    \nextgroupplot[ 
      ymin = 1e-11, ymax = 1e-2, ytickten = {-10, -8, ..., -2},
    ]
      \addplot[implicit residual] table [x="rank",y="res_impl"] {\adiDDDmsd};
      \addplot table [x="rank",y="res_expl"] {\adiDDDmsd};
      \addplot table [x="rank",y="res_expl"] {\adiSDDmsd};
      \addplot table [x="rank",y="res_expl"] {\adiSSDmsd};
      \addplot table [x="rank",y="res_expl"] {\adiSSSmsd};
    \nextgroupplot[ 
      ylabel = Time,
      ymin = 1e7, ymax = 1.5e11,
      ytickten = {7, ..., 12},
      yticklabels = {
        \qty{10}{\milli\second},
        \qty{100}{\milli\second},
        \qty{1}{\second},
        \qty{10}{\second},
        \qty{100}{\second},
      },
    ]
      \addplot table [x="rank",y="runtime_ns"] {\adiDDDrail};
      \addplot table [x="rank",y="runtime_ns"] {\adiSDDrail};
      \addplot table [x="rank",y="runtime_ns"] {\adiSSDrail};
      \addplot table [x="rank",y="runtime_ns"] {\adiSSSrail};
    \nextgroupplot[ 
      xlabel = Rank $g$,
      ymin = 1e7, ymax = 1.5e11,
      ytickten = {7, ..., 12},
    ]
      \addplot table [x="rank",y="runtime_ns"] {\adiDDDchip};
      \addplot table [x="rank",y="runtime_ns"] {\adiSDDchip};
      \addplot table [x="rank",y="runtime_ns"] {\adiSSDchip};
      \addplot table [x="rank",y="runtime_ns"] {\adiSSSchip};
    \nextgroupplot[ 
      ymin = 1e7, ymax = 1.5e11,
      ytickten = {7, ..., 12},
    ]
      \addplot table [x="rank",y="runtime_ns"] {\adiDDDmsd};
      \addplot table [x="rank",y="runtime_ns"] {\adiSDDmsd};
      \addplot table [x="rank",y="runtime_ns"] {\adiSSDmsd};
      \addplot table [x="rank",y="runtime_ns"] {\adiSSSmsd};
    \end{groupplot}
  \end{tikzpicture}

  \medskip
  \ref*{fig:num:random-rhs:legend} 
  \caption{%
    Implicit (top) and explicit residuals (middle)
    as well as wall-clock times (bottom)
    of \ac{ALE}~\eqref{eq:intro}
    over the rank~$g\in\N$ of a random constant term for
    ADI(\precZ, \precVR, \precTY) 
    as shown in \autoref{alg:adi};
    see \autoref{sec:num:random-rhs}.
    The middle row shows implicit residuals in the background~(%
    \protect\tikz[/pgfplots/every crossref picture] \protect\draw[implicit residual] plot coordinates { (0,0) };).
    The symbols S and D refer to IEEE single and double precision, respectively.
  }%
  \label{fig:num:random-rhs}
\end{figure*}

\subsection{Recover known solution}%
\label{sec:num:known-X}

In this experiment,
for a given rank~$\hat z\in\N$,
we first select random solution factors
$\hat Z\in\R^{n\times \hat z}$ and
$\hat Y\in\R^{\hat z\times \hat z}$ (symmetrized),
and set the right-hand side of \ac{ALE}~\eqref{eq:intro} to
a low-rank compression of
\begin{equation}
  \label{eq:rhs:known}
  G = \begin{bmatrix} E\hat Z & A\hat Z \end{bmatrix}
  \qquad\text{and}\qquad
  S = -\begin{bmatrix} & \hat Y \\ \hat Y \end{bmatrix}
\end{equation}
computed in IEEE double precision.
Observe that $\hat X := \hat Z \hat Y \hat Z\transpose$ indeed solves the resulting \ac{ALE}.
In contrast to \autoref{sec:num:gramian},
we do know that the solution can be exactly represented by a low-rank factorization,
moreover, using double precision.
Our main concern is not to solve an \ac{ALE} to the smallest possible residual,
but to approximate~$\hat X$ to the smallest possible error
\begin{equation}%
\label{eq:error}
  \errlabel.
\end{equation}
Therefore, the concerns raised in \autoref{rem:compression} do not apply,
and we will additionally report residuals for compressions of the \ac{ADI} approximations.

\Autoref{fig:num:known-x} shows the error in the approximation~$X$,
as well as the difference between the order of~$X$ to the known rank of~$\hat X$,
for varying ranks~$\hat z$.
Only the uniform double-precision \ac{ADI} is able to approximate the known solution to a small error,
uniformly over the ranks~$\hat z$ tested.
The error increases with the extent of single precision being used in the algorithm.
When accumulating the outer solution factor in single precision,
\mbox{ADI(S, D, D)},
the error is amplified by about one order of magnitude by low-rank compression.
Storing only the inner low-rank factors in double precision,
\mbox{ADI(S, S, D)},
shows no advantage over storing everything in single precision.

Recall that $\rank(X) \leq z$,
where $X=\ZYZ\transpose$,
such that \mbox{$z - \hat z$} describes the representation overhead of~$X$ w.r.t.~the known~$\hat X$.
\replaced{Some}{Most} approximations~$X$ for the Steel Profile (for \replaced{$\hat z\geq 40$}{$\hat z\geq 20$})
and Triple Chain examples (for $\hat z\geq 2$)
have almost full order ($z\approx n$ or even $z > n$).
Consequently,
storing the low-rank factors of $\ZYZ\transpose$ is more expensive than assembling the product.
The only exception is the uniform double-precision \ac{ADI} applied to the Triple Chain,
with $z-\hat z < 150$,
thus requiring significantly less storage after low-rank compression.

However, the Steel Profile and Triple Chain examples used are hardly large-scale.
Due to its larger dimension~$n=\num{20082}$,
the Chip example behaves rather forgiving:
all variants of the \ac{ADI} are truly low-rank ($z\leq\num{2700}\ll n$).
While low-rank compression does reduce the order~$z$,
it barely affects the error of the approximation.

In terms of runtime, yet again,
both \mbox{ADI($\ast$, D, D)} and both \mbox{ADI(S, S, $\ast$)} are the same.

\begin{figure*}
  \centering
  \begin{tikzpicture}
    \pgfplotstableread[col sep=comma]{Chip_0.0_/random-sol/ADI-DDD_compression=false_compression_interval=10_ignore_initial_guess=false_maxiters=50_reltol=1e-10_.csv}\adiDDDchip
    \pgfplotstableread[col sep=comma]{Chip_0.0_/random-sol/ADI-SDD_compression=false_compression_interval=10_ignore_initial_guess=false_maxiters=50_reltol=1e-10_.csv}\adiSDDchip
    \pgfplotstableread[col sep=comma]{Chip_0.0_/random-sol/ADI-SSD_compression=false_compression_interval=10_ignore_initial_guess=false_maxiters=50_reltol=1e-10_.csv}\adiSSDchip
    \pgfplotstableread[col sep=comma]{Chip_0.0_/random-sol/ADI-SSS_compression=false_compression_interval=10_ignore_initial_guess=false_maxiters=50_reltol=1e-10_.csv}\adiSSSchip
    \pgfplotstableread[col sep=comma]{FenicsRail_5177_/random-sol/ADI-DDD_compression=false_compression_interval=10_ignore_initial_guess=false_maxiters=50_reltol=1e-10_.csv}\adiDDDrail
    \pgfplotstableread[col sep=comma]{FenicsRail_5177_/random-sol/ADI-SDD_compression=false_compression_interval=10_ignore_initial_guess=false_maxiters=50_reltol=1e-10_.csv}\adiSDDrail
    \pgfplotstableread[col sep=comma]{FenicsRail_5177_/random-sol/ADI-SSD_compression=false_compression_interval=10_ignore_initial_guess=false_maxiters=50_reltol=1e-10_.csv}\adiSSDrail
    \pgfplotstableread[col sep=comma]{FenicsRail_5177_/random-sol/ADI-SSS_compression=false_compression_interval=10_ignore_initial_guess=false_maxiters=50_reltol=1e-10_.csv}\adiSSSrail
    \pgfplotstableread[col sep=comma]{TripleChain_100_/random-sol/ADI-DDD_compression=false_compression_interval=10_ignore_initial_guess=false_maxiters=200_reltol=1e-10_.csv}\adiDDDmsd
    \pgfplotstableread[col sep=comma]{TripleChain_100_/random-sol/ADI-SDD_compression=false_compression_interval=10_ignore_initial_guess=false_maxiters=200_reltol=1e-10_.csv}\adiSDDmsd
    \pgfplotstableread[col sep=comma]{TripleChain_100_/random-sol/ADI-SSD_compression=false_compression_interval=10_ignore_initial_guess=false_maxiters=200_reltol=1e-10_.csv}\adiSSDmsd
    \pgfplotstableread[col sep=comma]{TripleChain_100_/random-sol/ADI-SSS_compression=false_compression_interval=10_ignore_initial_guess=false_maxiters=200_reltol=1e-10_.csv}\adiSSSmsd
    \begin{groupplot}[
      group style = {
        columns = 3,
        rows = 3,
        horizontal sep = 2ex,
        vertical sep = 2ex,
        x descriptions at = edge bottom,
        y descriptions at = edge left,
      },
      scale only axis,
      xmin = 1,
      xmax = 50,
      enlarge x limits = true,
      xmode = log,
      ymode = log,
      ymin = 1e-11,
      ymax = 1e-2,
      ytickten = {-10, -6, -2},
      enlarge y limits = true,
      /tikz/only marks,
      samples at = {1,2,3,4,5,6,7,8,9,10,20,30,40,50,60,70,80,90,100},
      background/.style = {
        forget plot,
        opacity = 0.2,
        mark size = 4pt,
      },
      dimension/.style = {
        forget plot,
        mark = none,
        sharp plot, 
        dashed,
        samples = 2,
      },
    ]
    \nextgroupplot[
      title = Steel Profile,
      legend to name = fig:num:known-x:legend,
      legend columns = -1,
      legend entries = {%
        {ADI(\double, \double, \double)},
        {ADI(\single, \double, \double)},
        {ADI(\single, \single, \double)},
        {ADI(\single, \single, \single)},
      },
      ylabel = \errlabel,
      ytickten = {-10, -8, ..., -2},
    ]
      \addplot+[background] table [x="rank",y="err_tr"] {\adiDDDrail};
      \addplot              table [x="rank",y="err"]    {\adiDDDrail};
      \addplot+[background] table [x="rank",y="err_tr"] {\adiSDDrail};
      \addplot              table [x="rank",y="err"]    {\adiSDDrail};
      \addplot+[background] table [x="rank",y="err_tr"] {\adiSSDrail};
      \addplot              table [x="rank",y="err"]    {\adiSSDrail};
      \addplot+[background] table [x="rank",y="err_tr"] {\adiSSSrail};
      \addplot              table [x="rank",y="err"]    {\adiSSSrail};
    \nextgroupplot[
      title = Chip,
      ytickten = {-10, -8, ..., -2},
    ]
      \addplot+[background] table [x="rank",y="err_tr"] {\adiDDDchip};
      \addplot              table [x="rank",y="err"]    {\adiDDDchip};
      \addplot+[background] table [x="rank",y="err_tr"] {\adiSDDchip};
      \addplot              table [x="rank",y="err"]    {\adiSDDchip};
      \addplot+[background] table [x="rank",y="err_tr"] {\adiSSDchip};
      \addplot              table [x="rank",y="err"]    {\adiSSDchip};
      \addplot+[background] table [x="rank",y="err_tr"] {\adiSSSchip};
      \addplot              table [x="rank",y="err"]    {\adiSSSchip};
    \nextgroupplot[
      title = Triple Chain,
      ytickten = {-10, -8, ..., -2},
    ]
      \addplot+[background] table [x="rank",y="err_tr"] {\adiDDDmsd};
      \addplot              table [x="rank",y="err"]    {\adiDDDmsd};
      \addplot+[background] table [x="rank",y="err_tr"] {\adiSDDmsd};
      \addplot              table [x="rank",y="err"]    {\adiSDDmsd};
      \addplot+[background] table [x="rank",y="err_tr"] {\adiSSDmsd};
      \addplot              table [x="rank",y="err"]    {\adiSSDmsd};
      \addplot+[background] table [x="rank",y="err_tr"] {\adiSSSmsd};
      \addplot              table [x="rank",y="err"]    {\adiSSSmsd};
    \nextgroupplot[ 
      ylabel = $z - \hat z$,
      ymin = 1e1, ymax = 1.5e4, ytickten = {1, ..., 5},
    ]
      \addplot[dimension] { 5177 };
      \addplot+[background] table [x="rank",y="delta_rank_tr"] {\adiDDDrail};
      \addplot              table [x="rank",y="delta_rank"]    {\adiDDDrail};
      \addplot+[background] table [x="rank",y="delta_rank_tr"] {\adiSDDrail};
      \addplot              table [x="rank",y="delta_rank"]    {\adiSDDrail};
      \addplot+[background] table [x="rank",y="delta_rank_tr"] {\adiSSDrail};
      \addplot              table [x="rank",y="delta_rank"]    {\adiSSDrail};
      \addplot+[background] table [x="rank",y="delta_rank_tr"] {\adiSSSrail};
      \addplot              table [x="rank",y="delta_rank"]    {\adiSSSrail};
    \nextgroupplot[ 
      ymin = 1e1, ymax = 1.5e4, ytickten = {1, ..., 5},
    ]
      \addplot[dimension] { 20082 };
      \addplot+[background] table [x="rank",y="delta_rank_tr"] {\adiDDDchip};
      \addplot              table [x="rank",y="delta_rank"]    {\adiDDDchip};
      \addplot+[background] table [x="rank",y="delta_rank_tr"] {\adiSDDchip};
      \addplot              table [x="rank",y="delta_rank"]    {\adiSDDchip};
      \addplot+[background] table [x="rank",y="delta_rank_tr"] {\adiSSDchip};
      \addplot              table [x="rank",y="delta_rank"]    {\adiSSDchip};
      \addplot+[background] table [x="rank",y="delta_rank_tr"] {\adiSSSchip};
      \addplot              table [x="rank",y="delta_rank"]    {\adiSSSchip};
    \nextgroupplot[ 
      ymin = 1e1, ymax = 1.5e4, ytickten = {1, ..., 5},
    ]
      \addplot[dimension] { 602 };
      \addplot+[background] table [x="rank",y="delta_rank_tr"] {\adiDDDmsd};
      \addplot              table [x="rank",y="delta_rank"]    {\adiDDDmsd};
      \addplot+[background] table [x="rank",y="delta_rank_tr"] {\adiSDDmsd};
      \addplot              table [x="rank",y="delta_rank"]    {\adiSDDmsd};
      \addplot+[background] table [x="rank",y="delta_rank_tr"] {\adiSSDmsd};
      \addplot              table [x="rank",y="delta_rank"]    {\adiSSDmsd};
      \addplot+[background] table [x="rank",y="delta_rank_tr"] {\adiSSSmsd};
      \addplot              table [x="rank",y="delta_rank"]    {\adiSSSmsd};
    \nextgroupplot[ 
      ylabel = Time,
      ymin = 1e7, ymax = 1.5e11,
      ytickten = {7, ..., 12},
      yticklabels = {
        \qty{10}{\milli\second},
        \qty{100}{\milli\second},
        \qty{1}{\second},
        \qty{10}{\second},
        \qty{100}{\second},
      },
    ]
      \addplot table [x="rank",y="runtime_ns"] {\adiDDDrail};
      \addplot table [x="rank",y="runtime_ns"] {\adiSDDrail};
      \addplot table [x="rank",y="runtime_ns"] {\adiSSDrail};
      \addplot table [x="rank",y="runtime_ns"] {\adiSSSrail};
    \nextgroupplot[ 
      xlabel = Rank $\hat z$,
      ymin = 1e7, ymax = 1.5e11,
      ytickten = {7, ..., 12},
    ]
      \addplot table [x="rank",y="runtime_ns"] {\adiDDDchip};
      \addplot table [x="rank",y="runtime_ns"] {\adiSDDchip};
      \addplot table [x="rank",y="runtime_ns"] {\adiSSDchip};
      \addplot table [x="rank",y="runtime_ns"] {\adiSSSchip};
    \nextgroupplot[ 
      ymin = 1e7, ymax = 1.5e11,
      ytickten = {7, ..., 12},
    ]
      \addplot table [x="rank",y="runtime_ns"] {\adiDDDmsd};
      \addplot table [x="rank",y="runtime_ns"] {\adiSDDmsd};
      \addplot table [x="rank",y="runtime_ns"] {\adiSSDmsd};
      \addplot table [x="rank",y="runtime_ns"] {\adiSSSmsd};
    \end{groupplot}
  \end{tikzpicture}

  \medskip
  \ref*{fig:num:known-x:legend} 
  \caption{%
    Errors (top), representation overhead (middle), and wall-clock times (bottom)
    over the rank~$\hat z\in\N$ of a known solution~$\smash{\hat X}$ for
    ADI(\precZ, \precVR, \precTY) 
    as shown in \autoref{alg:adi};
    see \autoref{sec:num:known-X}.
    The symbols~\single{} and~\double{} refer to IEEE single and double precision, respectively.
    If visible, the dashed line denotes the system dimension~$n\in\N$.
    The solid markers denote uncompressed ADI approximates~$X$,
    opaque markers denote low-rank compressions of~$X$.
  }%
  \label{fig:num:known-x}
\end{figure*}

\section{Conclusions and future work}

We have investigated the low-rank \ac{ADI} applied to the algebraic Lyapunov equation.
We empirically observed that simply performing certain parts of the \ac{ADI} in a lower precision does not lead to highly accurate results.
We now revisit our hypotheses listed on \autopageref{list:hypotheses}:

\begin{enumerate}
  \item
    Our first hypothesis is partially true.
    The implicitly evaluated residual~\eqref{eq:residuals} suggests no deterioration of convergence
    \added{up to a target tolerance of about $\tau=10^{-8}$}
    even for the uniformly single-precision ADI\@.
    \added{%
    For smaller~$\tau$, solving the inner linear systems in single precision
    may lead to stagnation of the implicit residual 
    or vanishing of the \ac{ADI} increments.
    }\footnote{\added{%
    We abort the \ac{ADI} if $V_k=0$.
    Otherwise, this would cause the implicit residual to stagnate as well.
    }}%
    \deleted{All variants of the ADI terminate after the exact same number of iterations.}
    However, the corresponding explicit residuals are, in parts, substantially larger.

  \item
    Our second hypothesis is \replaced{wrong.}{true in terms of the residuals,
    but only excluding the Triple Chain example,
    and false w.r.t.~the wall-clock time.}
    \added{%
    The \mbox{ADI(S, D, D)} solution is several orders worse than the uniform double-precision one.
    }

    \added{In terms of wall-clock time,}
    accumulating the solution factor~$Z$ in single precision is,
    at least for the current implementation,
    not faster than accumulating the solution factor in double precision.
    This disadvantage is partly due to our current implementation,
    but also because we do not compress the low-rank solution in the end;
    see \autoref{rem:compression}.
    We lazily collect the individual low-rank factors~$V_k$ without concatenating them,
    but convert them to the desired precision.
    If the precisions \precZ{} and \precVR{} do not match,
    this requires an additional memory sweep,
    which can not be compensated for by a faster orthonormalization of~$Z$.

    However, the single-precision solution obtained from \mbox{ADI(S, D, D)} has the lowest residuals of all mixed- and single-precision variants of the \ac{ADI}.
    \deleted{%
    For all first-order examples tested,
    the residual is even comparable to the double-precision \ac{ADI}.
    }

  \item
    Our third hypothesis is true in terms of wall-clock time,
    but wrong w.r.t.~the residuals.
    Without low-rank compression,
    \mbox{ADI(S, S, D)} yields no better solutions than the uniform single-precision variant.
\end{enumerate}

\added{%
We would like to emphasize that,
excluding the Triple Chain example in \autoref{sec:num:random-rhs},
the \mbox{ADI(S, D, D)} solution is about 3 orders more accurate
than the single-precision solution,
while consuming the same amount of storage,
which is half the storage of the double-precision solution.
Furthermore,
the uniform single-precision \ac{ADI} can be a viable option, e.g.,
if one is only interested in derived quantities like the \Htwo{} system norm,
where one does not need to assemble the full solution.
In such cases, the single-precision \ac{ADI} provides better approximations than the single-precision Bartels-Stewart algorithm,
all at a fraction of the runtime.
}

An obvious further research direction is to improve the robustness of the \ac{ADI}.
To detect the stagnation of the explicitly evaluated residual~\eqref{eq:residuals},
even for the uniform double-precision \ac{ADI},
we recommend to evaluate the explicit residual every couple of \ac{ADI} iterations,
although algorithmically not necessary.
\added{Should implicit and explicit residual disagree too much, the \ac{ADI} can be aborted.}

\begin{remark}
  \added{%
  Lyapack~\cite{Pen00a} had a similar stagnation detection based on the explicit residual,
  in part because the implicit formulation of the residual~\cite{morBenKS13, morWolP16} was not known at the time.
  Since the implicit formulation is significantly cheaper to evaluate,
  the explicit formulation has seen little to no use in recent \ac{ADI} implementations.
  Instead, a more practical alternative could be to monitor the norm of the increment~$V_k T V_k\transpose$
  (and estimate the condition of the Lyapunov operator, as a byproduct of the shift computations)
  instead of the explicit residual.
  }
\end{remark}

Furthermore, the increment factor~$V_k$ could be scaled by some diagonal matrix,
or orthonormalized w.r.t.~all existing columns~$Z_{k-1}$,
before truncating to \precZ{}.
Either option would lose the Kronecker structure~\eqref{eq:kronecker} in the inner solution factor~$Y$.
The perhaps simplest option is, however,
to use a direct inner solver\footnote{ %
  Julia's~\cite{Julia-2017} built-in ``backslash'' resorts to an LU decomposition,
  which does not directly support single precision;
  the curious reader may check \lstinline!@edit lu(spzeros(Float32, 3, 3))!
  from within the Julia~REPL\@.
}
like MUMPS~\cite{AmeDL00}.\footnote{ %
  Thanks to
  Bastien Vieubl\'{e}\,\orcidlink{0000-0001-8429-7400} 
  for recommending MUMPS.jl~\cite{MUMPS.jl},
  which does support single-precision arithmetic.
}
\added{%
Another option, as is typical for mixed-precision applications,
is to add a self-correction mechanism like iterative refinement.
In the context of this paper,
iterative refinement is equivalent to restarting the \ac{ADI} with an initial guess.
Recall from \autoref{tab:gramian} that the single-precision \ac{ADI} can be $4\times$ as fast as the double-precision \ac{ADI}.
In light of the aforementioned stagnation detection,
as well as the current shift away from double precision in hardware development,
this gap is set to grow.
Thus, there certainly is headroom for the increased runtime cost incurred by higher-rank constant terms in the \ac{ALE}~\eqref{eq:intro} in subsequent \ac{ADI} runs~\cite[Remark~4.6]{SchS25}.
}

Further investigation is needed to be able to compress single- or multi-precision low-rank factorizations while not increasing their
explicit Lyapunov residuals~\eqref{eq:residuals}
or reconstruction error~\eqref{eq:error};
see \autoref{rem:compression}.
Besides exploring other methods computing the QR decomposition,
one should also explore other decompositions $Z=QF$,
where $Q$ is orthonormal as $F$ is square.
Of particular interest are multi-precision methods with $\precZ = \precision{Q} > \precision{F}$,
e.g.,~methods returning single-precision~$Q$ and double-precision~$F$.

A more general research direction is to revisit the finite-precision ADI as an inexact ADI~\cite{KueF20}.
For one, this field gives access to a more appropriate tolerance for the inner solver;
see \autoref{rem:GMRES:precision}.
More interestingly, however,
it should be possible to create an adaptive-precision \ac{ADI}
that switches from \mbox{ADI(S, D, D)} to \mbox{ADI(S, S, D)}.
That is, for later iterations,
the linear systems in \autoref{alg:adi:linsolve} of \autoref{alg:adi} can be solved in a lower precision.

Lastly, one could consider more than the three precisions,
e.g.,~$\precision{V} > \precision{R}$ instead of \precVR{}.
This either requires an implementation of GMRES
that is resilient to multi-precision input (see \autoref{rem:GMRES}),
or requires the \ac{ADI} to hold copies of matrix~$E$
in both precisions, \precision{V} and \precision{R}.

\addcontentsline{toc}{section}{Code and Data Availability}
\section*{Code and Data Availability}

\autoref{alg:adi} is implemented in Julia~\cite{Julia-2017}
and is available in DifferentialRiccatiEquations.jl~\cite{DRE.jl}.
The scripts specific to this research project and the data underlying the diagrams are available at:

\begin{center}
  \replaced{%
  \url{https://doi.org/10.5281/zenodo.18508895}
  }{%
  \url{https://doi.org/10.5281/zenodo.15808083}
  }
\end{center}

%% file: adi_algorithm.tex
\newcommand{\Zfact}[1]{\simplebox{fill1}{$Z_{#1}$}}
\newcommand{\Rfact}[1]{\simplebox{fill2}{$R_{#1}$}}
\newcommand{\Vfact}[1]{\simplebox{fill2}{$V_{#1}$}}
\newcommand{\Tfact}{\simplebox{fill3}{$T$}}
\newcommand{\Yfact}[1]{\simplebox{fill3}{$Y_{#1}$}}
\begin{algorithm*}
  \caption{%
    Mixed-Precision Low-Rank
    ADI(\simplebox{fill1}{\precZ}, \simplebox{fill2}{\precVR}, \simplebox{fill3}{\precTY}). 
    Outer solution factors are handled with precision~\simplebox{fill1}{\precZ},
    outer increment and residual factors with~\simplebox{fill2}{\precVR},
    inner low-rank factors with~\simplebox{fill3}{\precTY},
    where \begin{math}
      \simplebox{fill1}{\precZ} \geq
      \simplebox{fill2}{\precVR} \geq
      \simplebox{fill3}{\precTY}
    \end{math}.
  }%
  \label{alg:adi}
  \KwIn{%
    system matrices $A$, $E$, $G$, and $S$,
    parameters~$\{\alpha_0, \alpha_1, \ldots\}$,
    relative tolerance $\tau$
  }
  \KwOut{%
    $Z:=Z_{k+1}\in\C^{n\times z}$ and $Y:=Y_{k+1}\in\C^{z\times z}$
    such that $X \approx \ZYZ\transpose$ solves \ALE 
  }
  Initialize solution and residual factors:\linebreak
  \begin{math}
    \Zfact{0} \gets [\,],
    \qquad
    \Yfact{0} \gets [\,],
    \qquad
    \Rfact{0} \gets G,
    \qquad
    \Tfact \gets S
  \end{math}
  \;
  \For{$k \gets 0, 1, \ldots$}{%
    Compute increment factor
    \begin{math}
      \Vfact{k} \gets
      {(A\transpose +\alpha_{k} E\transpose)}^{-1}
      \Rfact{k}
    \end{math}
    \;\label{alg:adi:linsolve}
    Update residual factor
    \begin{math}
      \Rfact{k+1} \gets
      \Rfact{k}
      - 2\Re(\alpha_{k}) E\transpose
      \Vfact{k}
    \end{math}
    \;
    Augment solution factors:\linebreak
    \begin{math}
      \Zfact{k+1}
      \gets \Big[ \Zfact{k} \Big| \Vfact{k} \Big],
      \qquad
      \Yfact{k+1}
      \gets \begin{bmatrix}
        \Yfact{k} \\
        & -2\Re(\alpha_k)\Tfact
      \end{bmatrix}
    \end{math}
    \;
    \lIf{
      \begin{math}
        \Fnorm{
          \Rfact{k+1}
          \simplebox{fill3}{$T\vphantom{R_{k+1}}$}
          \Rfact{k+1}\transpose
        }
        \leq \tau \Fnorm{GSG\transpose}
    \end{math}
    }{break}
  }
\end{algorithm*}